\documentclass[9pt]{IEEEtran}

\usepackage{setspace}
\usepackage{array}
\usepackage{color}
\usepackage{amsmath,amssymb,mathrsfs}
\usepackage[mathscr]{eucal}
\usepackage{cite}
\usepackage{epsfig}
\usepackage{verbatim}
\usepackage{fixmath}
\usepackage{algorithm,algorithmic}
\usepackage{mathdots,xcolor}
\usepackage{bm,upgreek,nicefrac}
\usepackage{tikz}

\newtheorem{thm}{Theorem}[section]
\newtheorem{cor}[thm]{Corollary}
\newtheorem{lem}[thm]{Lemma}

\newtheorem{exmp}{Example}[section]

\newtheorem{rem}{Remark}[section]
\newcommand{\bi}[1]{\emph{\textbf{#1}}}
\newcommand{\LL} {\tilde{\mathcal{L}}_{\tilde{A}}^{(k)}(\lambda)}

\usepackage{cite}
\begin{document}
%
\title{Laplacian Controllability of Threshold Graphs}
%
%
%

\author{Shun-Pin~Hsu,~\IEEEmembership{Member,~IEEE,}
\thanks{Shun-Pin Hsu is with the Department
of Electrical Engineering, National Chung Hsing University, Taichung,
402, Taiwan, e-mail: (shsu@nchu.edu.tw).}}

\maketitle

\begin{abstract}
This paper is concerned with the controllability problem of a connected threshold graph following the Laplacian dynamics. An algorithm is proposed to generate a spanning set of orthogonal Laplacian eigenvectors of the graph from a straightforward computation on its Laplacian matrix. A necessary and sufficient condition for the graph to be Laplacian controllable is then proposed. The condition suggests that the minimum number of controllers to render a connected threshold graph controllable is the maximum multiplicity of entries in the conjugate of the degree sequence determining the graph, and this minimum can be achieved by a binary control matrix. The second part of the work is the introduction of a novel class of single-input controllable graphs, which is constructed by connecting two antiregular graphs with almost the same size. This new connecting structure reduces the sum of the maximum vertex degree and the diameter by almost one half, compared to other well-known single-input controllable graphs such as the path and the antiregular graph, and has potential applications in the design of controllable graphs subject to practical edge constraints. Examples are provided to illustrate our results.
\end{abstract}


%
\IEEEpeerreviewmaketitle

\section{Introduction}
%
%
%
%

Recently, many research efforts have been invested in the study of distributed and cooperative operations of multi-agent networks~\cite{lewis14,kno16,zhao17}. Special synchronization property and various consensus issues for heterogeneous or homogeneous systems with constraints or uncertainties were extensively studied~\cite{alv16,val17,alm17,wang17}. The simplicity of network units and scalability of network structure make such systems a promising solution for many engineering challenges in the automation, manufacturing, monitoring, intelligent transportation, and so on. Rather than on centralized coordination or control, a multi-agent system relies on the exchange of local information to reach certain system-wide goals. An immediate question is how the system leverages the network connectivity to ensure effective and efficient information propagation and complete assigned jobs? To set up a platform for the analysis, a frequently utilized measure is to consider the network from the perspective of a connected graph~\cite{mes10}, namely, using the nodes (or called vertices) to represent agents and edges to describe the local interactions of the network. Under this framework, the system evolution can be defined based on the specified dynamics. A fundamental issue that draws many attentions of researchers in the control theory community is the system controllability. This property directly reflects the effectiveness of control mechanism and is particularly essential in a large-scale system~\cite{liu11}. Suppose an agent is represented by a state variable, and the status of the agent, or the value of the variable, is directly influenced by its neighboring agents, namely, those agents directly connected to it. Under mild conditions such as the time-invariance of the connection parameters and the so-called consensus policy, the formulation leads naturally to the classical liner and time-invariant (LTI) control system evolving according to the Laplacian dynamics. In terms of multi-agent systems, this formulation is often called the leader-follower control dynamics. The leader agents are actually the independent input signals and maneuver the follower agents to particular status. While several methods including the well-known Kalman test are available, checking the Laplacian controllability of multi-agent systems defined on connected graphs is never trivial. The major challenges come from the numerical instability in dealing with the large number of state variables in complex systems. To overcome, graph-theoretic techniques are proposed to analyze the connection patterns and to identify the structure symmetry known as the \emph{equitable partition}. Existence of such symmetry serves as a sufficient condition for system uncontrollability~\cite{rah09}. A weaker version of the sufficient condition based on the so-called \emph{relaxed equitable partition}, or \emph{almost equitable partition} (AEP)~\cite{car07} was later proposed~\cite{cao13}. However, failing to satisfy these conditions leaves the Laplacian controllability of the connected graphs inconclusive. Following the line of the partition scheme, upper and lower bounds that make elegant use of the information on the AEP and on the so-called \emph{distance partition} respectively were presented for the dimension of controllable subspace of state variables~\cite{zha14}. These bounds are useful in the design of a Laplacian controllable graph. Unfortunately, the gap between the upper and lower bounds is narrow only in some extreme cases such as the \emph{path graph}. The process of exploring graph-theoretic methods based on partition schemes is still far from completed.

In practice, system architects might need to design various connection structures for networks.  To them, a sufficient condition for system controllability is preferable to that for uncontrollability. Since this preferred condition easily applicable to a large-sized graph is not available in general, special classes of graphs were studied and their controllability issues were partially or totally solved. Examples include the paths~\cite{par12}, multi-chain~\cite{hsu17a}, grids~\cite{not13}, circulant graphs~\cite{nab13}, and complete graphs~\cite{zhang11}. The eigenspaces of these graphs follow specific patterns and thus are analyzable using the well-known Popov-Belevitch-Hautus (PBH) test. Nevertheless, the class of known Laplacian controllable graphs are still insufficient and are in great need of expansion. In this work we continue to explore the class of controllable graphs and focus on threshold graphs, which have many applications in synchronizing process, cyclic scheduling, Guttman scales and so on~\cite{mahadev95}. We only consider the connected graphs throughout the paper to avoid the apparent uncontrollability due to the disconnection. A connected threshold graph can be uniquely determined from its degree sequence, or from a sequence of vertex add-ins via union or join operation starting from an isolated vertex. Its Laplacian eigenvalues, all integers, are actually the entries in the conjugate of its degree sequence~\cite{mer94}. Using this property, it was shown that for a threshold graph, it is single-input controllable if and only if it is an connected antiregular graph, namely, it is connected and has exactly two vertices with the same degree. The necessary and sufficient condition for a binary control vector to render the system controllable was also provided~\cite{agu15b}. This result was later extended to the multi-input case where more vertices with the same degree were considered~\cite{hsu16}. In this work we generalize these results to any connected and unweighted threshold graphs and study their controllability conditions under the general, binary and terminal control matrices respectively (a control matrix is terminal if each column of the matrix has only one nonzero entry). Moreover, we study a novel structure closely related to an antiregular graph. The study is motivated by the desire to find a single-input controllable graph that has different connecting statistics from that of the path or the antiregular graph~\cite{hsu17b}. In practice, to prevent a long response time, an upper bound might be imposed on the graph diameter; to reduce the operation or maintenance complexity of a single unit, an upper bound might be imposed on each vertex degree as well. A $k$-vertex path has its edge degrees bounded above by $2$, and its diameter by $k-1$. A $k$-vertex antiregular graph has the maximum edge degree $k-1$ and diameter $2$. In both cases the sum of these two parameters is $k+1$. We propose to combine two antiregular graphs with nearly equal number of vertices and reduce the sum to $\frac{k}{2}+4$ and $\frac{k+1}{2}+4$ for even and odd $k$ respectively, while ensuring the single-input controllability of the resulting graph.

The contributions of our results are twofold. Firstly, we propose a simple algorithm to generate a complete set of orthogonal Laplacian eigenvectors of a connected threshold graph. These vectors are obtained from a straightforward computation on the entries of the Laplacian matrix. Based on this result, we propose the necessary and sufficient condition for a connected threshold graph to be Laplacian controllable. In particular, we show that if the control matrix is binary, checking the controllability of a connected threshold graph is equivalent to checking the controllability of several intermediate connected threshold graphs generated in the entire construction process of the graph. Our result suggests that the difficulty in checking the controllability is not directly related to the number of vertices of the entire graph but mainly to that of equal degrees. The knowledge of Laplacian eigenspaces of the graph allows one to check the controllability locally and efficiently. Furthermore, we prove that under general and binary control matrices, the minimum number of controllers rendering the graph controllable is the same as the largest multiplicity of Laplcian eigenvalues; however, if the control matrix is terminal, the minimum number is the difference between the number of vertices and the number of vertices with different degrees. The second contribution of the paper is the introduction of a new class of single-input controllable graphs obtained by connecting two antiregular graphs. This novel class of graphs is unique in the following two senses: 1) it does not belong to any known family that has completely tractable Laplacian eigenspaces; 2) it has a much smaller sum of diameter and the maximum vertex degree, compared to that of a path or an antiregular graph. Our result demonstrates how adding one edge affects the Laplacian controllability of a graph, and thus has potential applications to the design of controllable graphs subject to edge constraints.

\section{PRELIMINARIES}
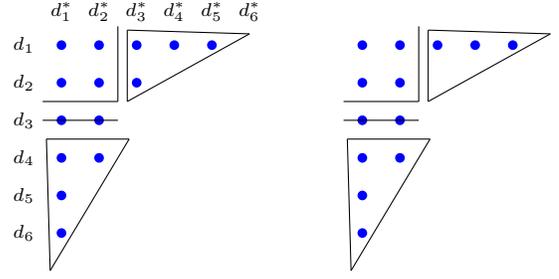
\begin{figure}
\begin{center}
	\begin{tikzpicture}
	\newcommand{\la}{0.5}	
	\newcommand{\lh}{0.5}
	\draw [blue,fill][thick] (0,0) circle [radius=0.05] node [black,below=4] {}; 
	\node[black] at (0,0.9*\lh) {\scriptsize $d_1^*$};
		\node[black] at (-\la,0) {\scriptsize $d_1$};
	\draw [blue,fill][thick] (\la,0) circle [radius=0.05] node [black,below=4] {}; 
	\node[black] at (\la,0.9*\lh) {\scriptsize $d_2^*$};
	\draw [blue,fill][thick] (2*\la,0) circle [radius=0.05] node [black,below=4] {}; 
		\node[black] at (2*\la,0.9*\lh) {\scriptsize $d_3^*$};
	\draw [blue,fill][thick] (2*\la,-\lh) circle [radius=0.05] node [black,below=4] {}; 
	\draw [blue,fill][thick] (\la,-3*\lh) circle [radius=0.05] node [black,below=4] {}; 
	\draw [blue,fill][thick] (3*\la,0) circle [radius=0.05] node [black,below=4] {}; 
		\node[black] at (3*\la,0.9*\lh) {\scriptsize $d_4^*$};
	\draw [blue,fill][thick] (4*\la,0) circle [radius=0.05] node [black,below=4] {}; 
		\node[black] at (4*\la,0.9*\lh) {\scriptsize $d_5^*$};
				\node[black] at (5*\la,0.9*\lh) {\scriptsize $d_6^*$};
	\draw [blue,fill] [thick](0,-\lh) circle [radius=0.05] node [black,below=4] {}; 
			\node[black] at (-\la,-\lh) {\scriptsize $d_2$};
	\draw [blue,fill] [thick](0,-2*\lh) circle [radius=0.05] node [black,below=4] {}; 
				\node[black] at (-\la,-2*\lh) {\scriptsize $d_3$};
	\draw [blue,fill] [thick](0,-3*\lh) circle [radius=0.05] node [black,below=4] {}; 
				\node[black] at (-\la,-3*\lh) {\scriptsize $d_4$};
	\draw [blue,fill] [thick](0,-4*\lh) circle [radius=0.05] node [black,below=4] {}; 
				\node[black] at (-\la,-4*\lh) {\scriptsize $d_5$};
	\draw [blue,fill] [thick](0,-5*\lh) circle [radius=0.05] node [black,below=4] {}; 
				\node[black] at (-\la,-5*\lh) {\scriptsize $d_6$};
	\draw [blue,fill][thick] (\la,-\lh) circle [radius=0.05] node [black,below=4] {}; 
	\draw [blue,fill] [thick](\la,-2*\lh) circle [radius=0.05] node [black,below=4] {}; 

	\draw [blue,fill][thick] (8*\la,0) circle [radius=0.05] node [black,below=4] {}; 
	\draw [blue,fill][thick] (9*\la,0) circle [radius=0.05] node [black,below=4] {}; 
	\draw [blue,fill][thick] (10*\la,0) circle [radius=0.05] node [black,below=4] {}; 
	\draw [blue,fill][thick] (9*\la,-3*\lh) circle [radius=0.05] node [black,below=4] {}; 
	\draw [blue,fill][thick] (11*\la,0) circle [radius=0.05] node [black,below=4] {}; 
	\draw [blue,fill][thick] (12*\la,0) circle [radius=0.05] node [black,below=4] {}; 
	\draw [blue,fill] [thick](8*\la,-\lh) circle [radius=0.05] node [black,below=4] {}; 
	\draw [blue,fill] [thick](8*\la,-2*\lh) circle [radius=0.05] node [black,below=4] {}; 
	\draw [blue,fill] [thick](8*\la,-3*\lh) circle [radius=0.05] node [black,below=4] {}; 
	\draw [blue,fill] [thick](8*\la,-4*\lh) circle [radius=0.05] node [black,below=4] {}; 
	\draw [blue,fill] [thick](8*\la,-5*\lh) circle [radius=0.05] node [black,below=4] {}; 
	\draw [blue,fill][thick] (9*\la,-\lh) circle [radius=0.05] node [black,below=4] {}; 
	\draw [blue,fill] [thick](9*\la,-2*\lh) circle [radius=0.05] node [black,below=4] {}; 

	\draw [black] [thin](\la+\la/2,-\lh-\lh/2) -- (-0.5*\la,-\lh-\lh/2);
	\draw [black] [thin](\la+0.5*\la,-\lh-0.5*\lh) -- (1.5*\la,0.5*\lh);
	\draw [black] [thin](\la+0.75*\la,-\lh-0.5*\lh) -- (1.75*\la,0.4*\lh);
	\draw [black] [thin](\la+4*\la,0.3*\lh) -- (1.75*\la,0.4*\lh);
	\draw [black] [thin](\la+4*\la,0.3*\lh) -- (1.75*\la,-1.5*\lh);
	\draw [black] [thin](\la+0.5*\la,-\lh*2) -- (-0.5*\la,-\lh*2);
	\draw [black] [thin](\la+0.8*\la,-\lh*2.5) -- (-0.4*\la,-\lh*2.5);
		\draw [black] [thin](-0.3*\la,-\lh*6) -- (-0.4*\la,-\lh*2.5);
				\draw [black] [thin](-0.3*\la,-\lh*6) -- (1.8*\la,-\lh*2.5);

	\draw [black] [thin](9*\la+\la/2,-\lh-\lh/2) -- (7.5*\la,-\lh-\lh/2);
	\draw [black] [thin](9*\la+0.5*\la,-\lh-0.5*\lh) -- (9.5*\la,0.5*\lh);
	\draw [black] [thin](9*\la+0.75*\la,-\lh-0.5*\lh) -- (9.75*\la,0.4*\lh);
	\draw [black] [thin](9*\la+4*\la,0.3*\lh) -- (9.75*\la,0.4*\lh);
	\draw [black] [thin](9*\la+4*\la,0.3*\lh) -- (9.75*\la,-1.5*\lh);
	\draw [black] [thin](9*\la+0.5*\la,-\lh*2) -- (7.5*\la,-\lh*2);
	\draw [black] [thin](9*\la+0.8*\la,-\lh*2.5) -- (7.6*\la,-\lh*2.5);
	\draw [black] [thin](7.7*\la,-\lh*6) -- (7.6*\la,-\lh*2.5);
	\draw [black] [thin](7.7*\la,-\lh*6) -- (9.8*\la,-\lh*2.5);
	\end{tikzpicture}
	\end{center}
\caption{The Ferrers-Sylverster diagrams for the degree sequences $\bi{d}_1=\{5,3,2,2,1,1\}$ (left) and $\bi{d}_2=\{5,2,2,2,1,1\}$ (right); the traces for both cases are $2$; the left diagram is symmetric in the sense that its bottom triangle can be flipped upward, with respect to the line connecting the two dots of $d_3$, and rotated $90^0$ clockwise to match its upper-right triangle, thus $\bi{d}_1$ determines a threshold graph; the right diagram does not have such a match thus $\bi{d}_2$ does not determine a threshold graph.}\label{fig:sylv}
\end{figure}

We define several standard notations and review some concepts used throughout the paper before starting our discussion. Let $\mathbb{R}$ and $\mathbb{N}$ be the sets of real and natural numbers, respectively, and $\textbf{1}$ and $\textbf{0}$  the column vectors of 1's and 0's, respectively. We use $I$ to represent the identity matrix with appropriate sizes and $\bi{e}_i$, $\bi{e}_{-i}$ are the $i$th and the $i$th to last column vectors, respectively, of $I$. Define the index set $\mathbb{I}_k:=\{1,2,\cdots k\}$. A vector $\bi{u}$ is called the \emph{$i$ to $j$} subvector of $\bi{v}=[\,v_1\,v_2\,\cdots\,v_n\,]^T\in\mathbb{R}^n$ if $\bi{u}=[\,v_i\,v_{i+1}\,\cdots\,v_j\,]^T$. For two sets $S_1$ and $S_2$, the set difference $S_1\setminus S_2$ is defined as $\{s|s\in S_1, s\notin S_2\}$. The floor and ceiling function of $x$, i.e., $\lfloor x\rfloor$ and $\lceil x\rceil$, are the largest integer not greater than $x$ and the smallest integer not less than $x$, respectively. Suppose $\lambda\in\mathbb{R}$, $\bi{v}\in\mathbb{R}^n$, and $P$ is a matrix of order $n$, meaning that $P\in\mathbb{R}^{n\times n}$. We call $(\lambda,\bi{v})$ an eigenpair of $P$ if $P\bi{v}=\lambda\bi{v}$. If $V=\mathbb{I}_k$ and $E$ is a subset of $\{\,(v_1,v_2)\,|\,v_1,v_2\in V\}$, then $(V,E)$ describes a $k$-vertex graph, or a graph on $k$ vertices, where $V$ and $E$ are called the vertex set and edge set respectively, of the graph. A path between vertices $v_1$ and $v_2$ is a subset $\{(v_1,u_1),(u_1,u_2),(u_2,u_3),\cdots,(u_{m-1},u_m),(u_m,v_2)\}$ of $E$ where $u_1,u_2,\cdots,u_m\in V$. A graph is connected if for every pair of vertices there exists a path between them. A graph is unweighted if all its connecting edges carry the same weight. A graph is undirected if every entry of its edge set is not an ordered pair. A unweighted and undirected graph is simple if it has no (self) loops and no multiple edges~\cite{god01}. In this work we restrict our discussion only to the class of simple and connected graphs. In a $k$-vertex simple graph described by $(V,E)$, if $v_1,v_2\in V$ and $(v_1,v_2)\in E$, $v_1$ and $v_2$ are called neighbors. The neighbor set $\mathcal{N}_v$ of the vertex $v$ is $\{u\,|\,(v,u)\in E\}$. The degree of vertex $v$ is defined as the cardinality of $\mathcal{N}_v$, written as $|\mathcal{N}_v|$, meaning the number of elements in $\mathcal{N}_v$. A vertex $v$ is called an isolated vertex, a terminal vertex, and a dominating vertex if $|\mathcal{N}_v|$ is $0,1$ and $k-1$, respectively. Furthermore, define $\mathcal{A},\mathcal{D},\mathcal{L}\in\mathbb{R}^{k\times k}$ where $\mathcal{D}$ is the diagonal matrix with its $i$th diagonal term being the degree of the $i$th vertex; $\mathcal{A}$ is the adjacent matrix of the graph, meaning that its $(i,j)$th element is $1$ if $(i,j)\in E$ and is $0$ otherwise; $\mathcal{L}:=\mathcal{D}-\mathcal{A}$ is the Laplacian matrix of the graph. Suppose $d_i$ is the degree of the $i$th vertex of a $k$-vertex graph and satisfies $d_i\geq d_{i+1}$ for each $i\in\mathbb{I}_{k-1}$, then the degree sequence of the graph is
\begin{equation}\label{eq:tdseq}
\begin{split}
\bi{d}:&=(d_1,d_2,\cdots,d_k)\\
&=\left(\tilde{d}_1,\cdots,\tilde{d}_1,\tilde{d}_2,\cdots,\tilde{d}_2,\cdots,\tilde{d}_{\tilde{k}},\cdots,\tilde{d}_{\tilde{k}}\right) \end{split}
\end{equation}
where $\tilde{d}_i$, with multiplicity $m_i$, satisfies $\tilde{d}_i >\tilde{d}_{i+1}$ for each $i$ in  $\mathbb{I}_{\tilde{k}-1}$. The conjugate of the degree sequence $\bi{d}$ is
\begin{equation}\label{eq:tdsseq}
\begin{split}
\bi{d}^*&:=(d_1^*,d_2^*,\cdots,d_k^*)\\
&=\left(\tilde{d}_1^*,\cdots,\tilde{d}_1^*,\tilde{d}_2^*,\cdots,\tilde{d}_2^*,\cdots,\tilde{d}_{\tilde{k}^*}^*,\cdots,\tilde{d}_{\tilde{k}^*}^*\right)
\end{split}
\end{equation}
where $d_i^*:=|\{\,j\,|\,d_j\geq i\,\}|$ and $\tilde{d}_i^*$, with multiplicity $m_i^*$, satisfies $\tilde{d}_i^*>\tilde{d}_{i+1}^*$ for each $i$ in $\mathbb{I}_{\tilde{k}^*-1}$. Also, we define $\bar{m}_0=\bar{m}_0^*:=0$ and for $i\ge 1$,
$\bar{m}_i:=\sum_{j=1}^im_j$, $\bar{m}_i^*:=\sum_{j=1}^im_j^*$. The trace of a degree sequence $\bi{d}$ is $\tau_{\bi{d}}:=|\{j\,|\,d_j\geq j\,\}|$. We say a sequence $\bi{d}$ of nonnegative integers is graphical if there exists a graph whose degree sequence is $\bi{d}$. Suppose  $\bi{d}=(d_1,d_2,\cdots,d_k)$. It was shown that the necessary and sufficient condition for $\bi{d}$ to be graphical is~\cite[p.72]{mol12}
\begin{equation}\label{lem:gra}
\sum_{i=1}^j(d_i+1)\leq\sum_{i=1}^jd_i^*, \quad\mbox{$\forall j\in\mathbb{I}_{\tau_\bi{d}}$}.
\end{equation}
In the extreme case that the equality in~(\ref{lem:gra}) holds, the degree sequence determines the so called \emph{threshold graphs} or \emph{maximal graphs}~\cite{mer94}. Threshold graphs admit several equivalent definitions~\cite[Theorem~1.2.4]{mahadev95}. One of the simplest ways to identify them is using the so-called Ferrers-Sylvester diagram. In~Fig.~\ref{fig:sylv} two degree sequences are checked using the diagram to see if they form threshold graphs.
Another well-known definition is based on recursive \emph{vertex add-in} operations to form the graph. Specifically, if we start from an isolated vertex and increase the number of vertices one by one, with one of the two following operations:
\begin{enumerate}
\item adding an isolated vertex to the graph (union);
\item adding a dominating vertex to the graph (join);
\end{enumerate}
the resulting graph is a threshold graph~\cite{ban17}. We use $\mathcal{C}_T$ to represent the class of connected threshold graphs. A $k$-vertex graph in $\mathcal{C}_T$ is written as $\mathbb{G}_T^{(k)}$, and its corresponding Laplacian matrix $\mathcal{L}_T^{(k)}$. The constructing process of $\mathbb{G}_T^{(k)}$ can be encoded with a binary string $\mathcal{S}^{(k)}$ of size $k-1$, whose $i$th element is $1$ if the $(i+1)$th vertex is added via the join operation, and is $0$ if it is via the union operation~\cite{bap13}. Let $\mathcal{S}^{(k)}=(\,s_1,s_2,\cdots,s_{k-1}\,)$. Then we can write $(\,s_1,s_1,s_2,s_3,\cdots,s_{k-1}\,)=
(\,\tilde{s}_1,\cdots,\tilde{s}_1,\tilde{s}_2,\cdots,\tilde{s}_2,\cdots, \tilde{s}_{\tilde{k}},\cdots,\tilde{s}_{\tilde{k}}\,)$ where $\tilde{s}\neq\tilde{s}_{i+1}$ and the multiplicity of $\tilde{s}_i$, written as $m_{\sigma_i}$, is a permutation of $m_i$ in~(\ref{eq:tdseq}) for $i\in\mathbb{I}_{\tilde{k}}$. Naturally, the $\mathcal{S}^{(k)}$ generating $\mathbb{G}_T^{(k)}$ ends with $1$. A succession $(\,s_i,s_{i+1},\cdots,s_j\,)$ of $s\in\{0,1\}$ in $\mathcal{S}^{(k)}=(\,s_1,s_2,\cdots,s_{k-1}\,)$ means that $s_i=\cdots=s_j=s$ and $s_{i-1},s_{j+1}$ either not defined or not equal to $s$; in addition, if $i=1$ then $j\ge i$, otherwise, $j>i$. $\mathcal{S}$ is called a substring of $\mathcal{S}^{(k)}$ if $\mathcal{S}$ generates $\mathbb{G}_T^{(|\mathcal{S}|+1)}$ whose vertices are those in the first $|\mathcal{S}|+1$ vertices of $\mathbb{G}_T^{(k)}$ generated by $\mathcal{S}^{(k)}$. Observe that $\mathcal{S}^{(k)}$ induces the substrings $\mathcal{S}^{(k_1)},\mathcal{S}^{(k_2)},\cdots$ that generate the connected threshold subgraphs $\mathbb{G}_T^{(k_1)},\mathbb{G}_T^{(k_2)},\cdots$ of $\mathbb{G}_T^{(k)}$, respectively, where $\mathcal{S}^{(k_1)}\subseteq\mathcal{S}^{(k_2)}\subseteq\cdots\subseteq\mathcal{S}^{(k)}$ for some $k_1\le k_2\le\cdots\le k$. We are particularly interested in the substrings $\mathcal{S}^{(k_i)}$, $i\in\mathbb{I}_{\bar{s}}$, where $\bar{s}$ is the number of successions in $\mathcal{S}^{(k)}$ and the terminal $1$ of each $\mathcal{S}^{(k_i)}$ is preceded by a succession of $0$, or is the final entry of a succession of $1$. Thus $\mathbb{G}^{(k_i)}_T$ is the smallest connected threshold subgraph of $\mathbb{G}^{(k)}_T$ with its constructing string $\mathcal{S}^{(k_i)}$ containing the first $i$ successions of $\mathcal{S}^{(k)}$. $\mathcal{S}^{(k_i)}$ with this property is called an essential substring of $\mathcal{S}^{(k)}$. A threshold graph has an integral Laplacian spectrum and tractable eigenspace since its construction involves only the graph operations of union and join~\cite{mer94,mer98}. These properties are intimately related to our main result of generating orthogonal Laplacian eigenvectors directly from the Laplacian matrix of the graph. A special threshold graph, which allows only interlacing operations of join and union in the constructing process, or equivalently, has exactly two vertices of the same degree, is known as the antiregular graph~\cite{mer03}. The class of connected antiregular graphs was proved to be the only class of threshold graphs controllable by single controllers~\cite{agu15b}. Let $\mathcal{C}_A$ be the class of these special graphs. A $k$-vertex graph in $\mathcal{C}_A$ is written as $\mathbb{G}_A^{(k)}$  and its Laplacian matrix $\mathcal{L}_A^{(k)}$.

Consider an autonomous linear and time-invariant (LTI) system defined on a $k$-vertex graph $(V,E)$, in the sense that each system states $x_i$ evolves according to the so-called consensus policy:
\begin{equation}\label{eq:cons}
\dot{x}_i=-\sum_{j\in\mathcal{N}_i}(x_i-x_j)
\end{equation}
for each vertex $i\in V=\mathbb{I}_k$, where $\mathcal{N}_i\subseteq V$, is the neighbor set of vertex $i$. In the matrix form we have
\begin{equation}\label{eq:automcons}
\dot{\bi{x}}=-\mathcal{L}\bi{x}
\end{equation}
where $\bi{x}:=[\,x_1\;x_2\;\cdots\;x_k\,]^T$ and $\mathcal{L}$ is the Laplacian matrix of the graph
This is known as the Laplacian dynamics~\cite[p.1613]{agu15} on the graph $(V,E)$. Suppose $p$ control inputs $\bi{u}(t)=\left[\,u_1(t)\;u_2(t)\;\cdots\; u_p(t)\,\right]^T$ is applied to~(\ref{eq:automcons}) via the $p$-control matrix $B:=[\,\bi{b}_1\,\bi{b}_2\,\cdots,\bi{b}_p\,]$ where the $(i,j)$th element of $B$, or $b_{ij}$, is
nonzero if vertex $i$ is connected to input $u_j(t)$, and is $0$ otherwise. $B$ is called binary if the nonzero entry is restricted to $1$, otherwise it is called general. In particular, $B$ is called terminal if every column has only one nonzero entry. The system under control evolves according to the following equation:
\begin{equation}\label{eq:mcons}
\dot{\bi{x}}=-\mathcal{L}\bi{x}+B\bi{u}(t),
\end{equation}
which is written as $(\mathcal{L},B)$ for simplicity. In the single-input case, namely $p=1$, the more specific notation $(\mathcal{L},\bi{b})$ is used. Clearly $\bi{b}_j$ records the connection information of controller $j$ to each vertex. We say $\hat{\bi{b}}_j$ is a local control vector of $\bi{b}_j$ subject to some vertices if $\hat{\bi{b}}_j$ records only the connection information of controller $j$ to these vertices. Let $\mathcal{S}$ be a substring of $\mathcal{S}^{(k)}$ that generates $\mathbb{G}^{(k)}_T$ with Laplacian matrix $\mathcal{L}_T^{(k)}$ that is under the $p$-control matrix $B$. If $\mathcal{S}$ generates $\mathbb{G}_T^{(|\mathcal{S}|+1)}$, then the local control matrix $\hat{B}(\mathcal{S})\in\mathbb{R}^{(|\mathcal{S}|+1)\times p}$ is the submatrix of $B$ where the $j$th column of $\hat{B}(\mathcal{S})$ is the local control vector $\hat{\bi{b}}_j$ subject to vertices of $\mathbb{G}_T^{(|\mathcal{S}|+1)}$. Suppose $\bi{v}$ is en eigenvector of $\mathcal{L}_T^{(k)}$. $\hat{\bi{v}}$ is called the subvector of $\bi{v}$ subject to $\mathbb{G}_T^{(|\mathcal{S}|+1)}$ if $\hat{\bi{v}}$ is an eigenvector of $\mathcal{L}_T^{(|\mathcal{S}|+1)}$. Moreover, we define the local checking matrix $C(\mathcal{S})\in\mathbb{R}^{(|\mathcal{S}|+1)\times m}$ where $m$ is the number of entries of last succession in $\mathcal{S}$. If the last succession of $\mathcal{S}$ is a succession of $1$, the $j$th column of $C(\mathcal{S})$ is $\bi{e}_1-\bi{e}_{j+1}$ for each $j\in\mathbb{I}_{m-1}$. The last column of $C(\mathcal{S})$ is $\bi{e}_1-\bi{e}_{m+1}$ if $\mathcal{S}$ has only succession and is $(|\mathcal{S}|+\delta-m)\bi{e}_1-\sum_{j=m+1}^{|\mathcal{S}|+\delta}\bi{e}_j$ where $\delta=1$ otherwise.
If the last succession of $\mathcal{S}$ is a succession of $0$, all $\bi{e}_i$'s are replaced with $\bi{e}_{-i}$'s and $\delta$ with $0$. In the next section we study the Laplacian controllability problems of connected threshold graphs. The study is based on the following conclusion from the classical Popov-Belevitch-Hautus (PBH) test.
 \begin{thm}\label{thm:ns}~\cite{chen99} The system $(\mathcal{L},B)$ defined in~(\ref{eq:mcons}) is uncontrollable if and only if $\mathcal{L}$ has an eigenvector in the null space of $B^T$.
\end{thm}

 \section{Main Results}
 \subsection{Laplacian Eigenspaces of Threshold Graphs}
Consider a $k$-vertex connected threshold graphs with degree sequence $\bi{d}=(d_1,d_2,\cdots,d_k)$. Note that a connected threshold graph can be determined uniquely from its degree sequence (up to isomorphism). Following the definition of the degree sequence and the symmetric property of the graph in the Ferrers-Sylvester diagram, we can write its Laplacian matrix $\mathcal{L}_T^{(k)}$, without loss of generality, as
\begin{equation}\label{def:lkt}
-\left[\begin{array}{ccccccc}
-d_1	&1	&\cdots	&1	&1	&\cdots  &1\\
1 &\ddots	&\ddots &\vdots	&x	&\cdots 	&x\\
\vdots &\ddots	&\ddots &1	&\vdots 	&\vdots &\vdots\\
1 &\cdots	&1 &-d_{\tau_{\bi{d}}+1}	&x	&\vdots 	&\vdots\\
1 &x	&\cdots &x  &-d_{\tau_{\bi{d}}+2}	&\ddots 	&\vdots\\
\vdots &\vdots	&\cdots &\cdots &\ddots &\ddots &x \\
 1 &x	&\cdots &\cdots &\cdots &x &-d_k
\end{array}\right]
\end{equation}
where $x$ in the $(i,j)$th entry is $1$ if vertices $i$ and $j$ are connected, otherwise it is $0$. If the $(i,j)$th element $\ell_{i,j}$ of $\mathcal{L}_T^{(k)}$ is $1$, then $i>j$ implies $\ell_{i,m}=1$ for each $m\in\mathbb{I}_{j-1}$, and $i<j$ implies
$\ell_{m,j}=1$ for each $m\in\mathbb{I}_{i-1}$. Similarly, if $\ell_{i,j}=0$, then $i>j$ implies $\ell_{m,j}=0$ for each $m\in\{i+1,i+2,\cdots,k\}$, and $i<j$ implies $\ell_{i,m}=0$ for each $m\in\{j+1,j+2,\cdots,k\}$. To see this we can use the definition of sequential~\emph{union or join} operations to construct the graph and manage to align the vertices, from left to right, when the construction is complete. The rule to decide the order of each vertex in the line is: the vertex is added to the leftmost and rightmost locations of the current line of vertices if it is added via the join and union operations respectively. When we complete the $k$-vertex graph, we obtain a line of vertices with their degrees, from left to right, being decreasing. It turns out that a vertex with degree $t$ where $t<k$, has those, excluding the vertex itself, with top $t$ degrees as its neighbors due to the sure join operations performed after the $t$-edge vertex is added and possible join operations when it is added. In the following we first review and present related results that help further analysis of the Laplacian eigenspace of a connected threshold graph.
\begin{thm}(cf.\cite{mer94})~\label{thm:dseq}
Let $\tau_{\bi{d}}$ be the trace of $\bi{d}$ with corresponding $\mathcal{L}_T^{(k)}$ in~(\ref{def:lkt}). Let $\bi{v}:=[\,v_1\,\cdots\,v_k\,]^T$ be the $(\tau_{\bi{d}}+1)$th column of $\mathcal{L}^{(k)}_T$, and for $i\in\mathbb{I}_k$, $m_i$ and $m_i^*$ are the multiplicities defined  in~(\ref{eq:tdseq}) and ~(\ref{eq:tdsseq}) respectively. We then have
\begin{enumerate}
\item \begin{equation}
d_i=\left\{\begin{array}{ll}
d_i^*-1 &\mbox{if $i\in\mathbb{I}_{\tau_{\bi{d}}}$}\\\tau_{\bi{d}} &\mbox{if $i=\tau_{\bi{d}}+1$}\\
d_{i-1}^* &\mbox{o.w.}
\end{array}\right.;
\end{equation}
\item $d_{\tau_{\bi{d}}}>d_{\tau_{\bi{d}}+2}$;
\item the multiplicity of $d_{\tau_{\bi{d}}+1}$ is at least $2$;
\item \begin{equation}
v_i=\left\{\begin{array}{rl}-1 &\mbox{if $i\in\mathbb{I}_{\tau_{\bi{d}}}$ }\\
\tau_{\bi{d}} &\mbox{if $i=\tau_{\bi{d}}+1$}\\
0 &\mbox{o.w.}.\end{array}\right.;
\end{equation}
\item \begin{equation*}
m_i^*=\left\{\begin{array}{ll}
m_i-1 &\mbox{if  $d_{\tau_{\bi{d}}+1}=\tilde{d}_i$}\\
m_i &\mbox{o.w.}
\end{array}\right..
\end{equation*}
\end{enumerate}
\end{thm}
\begin{IEEEproof}
In the first item, the case of $i\in\mathbb{I}_k\setminus\{\tau_{\bi{d}}+1\}$ follows directly from the symmetry of the Ferrers-Sylvester diagran of a threshold graph. If $d_{\tau_{\bi{d}}+1}>\tau_{\bi{d}}$, then the trace of the graph will be at least $\tau_{\bi{d}}+1$, a contradiction. If
$d_{\tau_{\bi{d}}+1}<\tau_{\bi{d}}$, $d_{\tau_{\bi{d}}}^*=\tau_{\bi{d}}$ and thus $d_{\tau_{\bi{d}}}=\tau_{\bi{d}}$ due to the symmetry of the diagram. On the other hand, the definition of a threshold graph requires that $d_i^*=d_i+1$  for each $i\in\mathbb{I}_{\tau_{\bi{d}}}$, which leads to a contradiction. In the second item, note that elements of $\bi{d}$ are in non-increasing order. If $d_{\tau_{\bi{d}}}=d_{\tau_{\bi{d}}+2}$, then the first part of the lemma implies $d_{\tau_{\bi{d}}}=d_{\tau_{\bi{d}}+1}=d_{\tau_{\bi{d}}+2}=\tau_{\bi{d}}$, which breaks the symmetry of the Ferrers-Sylvester diagram of a threshold graph. In the third item, if the multiplicity of $d_{\tau_{\bi{d}}+1}$ is $1$, then the first part of the lemma implies $d_{\tau_{\bi{d}}}\geq\tau_{\bi{d}}+1$ and $d_{\tau_{\bi{d}}+2}\leq\tau_{\bi{d}}-1$, which again breaks the symmetry of the Ferrers-Sylvester diagram of a threshold graph. Item 4) follows directly from the position property of $1$ and $0$ in~(\ref{def:lkt}) and the fact that each column sum of $\mathcal{L}^{(k)}_T$ is $0$. The final item is a natural result from the Ferrers-Sylvester diagram of a threshold graph.
\end{IEEEproof}

Now we present a systematic approach to generate a set of $k-1$ eigenvectors of $\mathcal{L}_T^{(k)}$. This approach includes two steps only. The first step is to modify the upper-right triangular part of $\mathcal{L}^{(k)}_T$. The second step is to make sure the column sum after the modification is $0$. Specifically,
\begin{enumerate}
\item Let $v_{ij}$ be the $(i,j)$th entry of $\mathcal{V}$, and let $\mathcal{V}$ be $\mathcal{L}^{(k)}_T$ initially. If $j>i$, then $v_{ij}$ is replaced with $-1-v_{ij}$. As a result, $v_{ij}$ becomes $0$ if it is originally $-1$; or becomes $-1$ if it is originally $0$.
\item Update the diagonal terms of $\mathcal{V}$ generated by the first step, such that each column sum is $0$.
\end{enumerate}
This approach is algorithmized in Algorithm~1. Note that in the algorithm $v_{ii}$ is updated \emph{before} $v_{st}$, where $t>s, t>i$, is replaced such that two nested \textbf{for} loops are enough. Observe that for generating the matrix $\mathcal{V}$ from $\mathcal{L}^{(k)}_T$, Algorithm~1 modifies: 1). ${k}\choose{2}$ entries, for each using only one \emph{addition} operation; 2). $k-1$ entries, for each using $k-2$ \emph{addition} operations for summation. Thus for a threshold graph with data size $n$, we can run the algorithm with complexity $\mathcal{O}(n)$. The following theory verifies that the $i$th column of the output $\mathcal{V}$ of Algorithm~1 is an eigenvector corresponding to the $i$th largest eigenvalue of the given $\mathcal{L}_A^{(k)}$. Along with the trivial eigenvector $\textbf{1}$, this set of $k$ eigenvectors can be easily verified to be orthogonal. This result is based on the well-known fact that adding a dominating vertex to a $(k-1)$-vertex simple and connected graph results in an eigenpair $(k,\bi{v})$ where $\bi{v}=[\,(k-1)\,-1\,\cdots\,-1\,]^T$ is the first column of the Laplacian matrix of the resulting $k$-vertex graph and the first diagonal entry of the matrix corresponds to the degree of the added dominating vertex~\cite{mer98}.

 \begin{algorithm}[t!]
 \caption{Algorithm for generating eigenvectors of $\mathcal{L}_T^{(k)}$}
 \begin{algorithmic}[1]
 \renewcommand{\algorithmicrequire}{\textbf{Input: }}
 \renewcommand{\algorithmicensure}{\textbf{Output:}}
 \REQUIRE $\mathcal{L}_T^{(k)}$, the Laplacian matrix corresponding a $k$-vertex threshold graph
 \ENSURE  $\mathcal{V}$, satisfying $\mathcal{L}_T^{(k)}\mathcal{V}=\mathcal{V}D^*$, where $D^*$ is diagonal with $d^*_i$, the $i$th component of degree conjugate $\bi{d}^*$, being its $(i,i)$th entry.
 \\ \textit{Initialisation} : Let $v_{ij}$ be the $(i,j)$th entry of $\mathcal{V}$.
  \STATE $\mathcal{V}\leftarrow\mathcal{L}_T^{(k)}$
  \FOR {$j=2$ to $k$}
  \FOR {$i=1$ to $j-1$}
  \STATE $v_{ij}\leftarrow -1-v_{ij}$
    \ENDFOR
    \STATE $v_{jj}\leftarrow -\sum_{i=1,i\neq j}^kv_{ij}$
    \ENDFOR
    \STATE Remove the zero column of $\mathcal{V}$
    \STATE Append $\textbf{1}$, the column of $1$'s, to the last column of $\mathcal{V}$
 \RETURN $\mathcal{V}$
 \end{algorithmic}
 \end{algorithm}

\begin{thm} \label{thm:evec} Suppose a threshold graph with its Laplacian matrix $\mathcal{L}_T^{(k)}$ in~(\ref{def:lkt}) has the degree sequence $\bi{d}=(d_1,\cdots,d_k)$ and its conjugate $\bi{d}^*=(d_1^*,\cdots,d_k^*)$. Then the $j$th column of $\mathcal{V}$ generated by Algorithm~1 is the eigenvector of $\mathcal{L}_T^{(k)}$ corresponding to $d_j^*$, the $j$th largest eigenvalue of $\mathcal{L}^{(k)}_T$.
\end{thm}
\begin{IEEEproof} Consider the case that the eigenvalue is in the set of $\{d_1^*,d_2^*,\cdots,d_{\tau_{\bi{d}}}^*\}$.
From the Ferrers-Sylvester diagram, it can be seen that if the vertex with the maximum degree is removed (along with its connecting edges), the remaining graph becomes a $(d_2^*-1)$-vertex threshold graph. Let
\begin{equation}
\mathcal{P}^{(m)}=\left[\begin{array}{rrrr}
d_1	&-1	&\cdots	&-1\\
-1 &1	 & & \\
\vdots & &\ddots &\\
-1 & &&1
\end{array}\right]_{m\times m},
\end{equation}
\begin{equation}
\mathcal{P}^{(m)}_r=\left[\begin{array}{rrr}
\textbf{O}_{(r-1)\times (r-1)}	&&\\
 &\mathcal{P}^{(m)}	 &\\
&&\textbf{O}
\end{array}\right]_{k\times k}
\end{equation}
and
\begin{equation}
\mathcal{Q}^{(m)}_r=\left[\begin{array}{rrr}
\textbf{O}_{(r-1)\times (r-1)}	&&\\
 &\mathcal{L}_T^{(m)}	 &\\
&&\textbf{O}
\end{array}\right]_{k\times k}.
\end{equation}
We can thus write
$\mathcal{L}_T^{(k)}=\mathcal{L}_T^{(d_1^*)}=\mathcal{P}_1^{(d_1^*)}+\mathcal{Q}_2^{(d_2^*-1)}$.
Let $\bi{v}_2$ be the second column of $\mathcal{Q}_2^{(d_2^*-1)}$. Then
$\bi{v}_2^T=[\,0\quad d_2^*-2\quad-1\,\,\cdots\,-1\quad\bi{0}^T\,]$.
Clearly, $\mathcal{Q}_2^{(d_2^*-1)}\bi{v}_2=(d_2^*-1)\bi{v}_2$. Observe that $\mathcal{P}_1^{(d_1^*)}\bi{v}_2=\bi{v}_2$. Thus $\mathcal{L}^{(k)}_T\bi{v}_2=(1+d_2^*-1)\bi{v}_2=d_2^*\bi{v}_2$.
Similarly, we can write $\mathcal{Q}_2^{(d2^*-1)}=\mathcal{P}_2^{(d_2^*-1)}+\mathcal{Q}_3^{(d_3^*-2)}$. Let $\bi{v}_3$ be the third column of $\mathcal{Q}_3^{(d_3^*-2)}$. Then $\mathcal{Q}_3^{(d_3^*-2)}\bi{v}_3=(d_3^*-2)\bi{v}_3$. Observe that
$\mathcal{P}_1^{(d_1^*)}\bi{v}_3=\bi{v}_3$ and $\mathcal{P}_2^{(d_2^*-1)}\bi{v}_3=\bi{v}_3$. Thus $\mathcal{L}^{(k)}_T\bi{v}_3=(2+d_3^*-2)\bi{v}_3=d_3^*\bi{v}_3$. Continuing this argument we conclude that if $\bi{v}_{\ell}$ is the $\ell$th column of $\mathcal{Q}_{\ell}^{(d_{\ell}^*-(\ell-1))}$ for each $\ell\in\mathbb{I}_{\tau_{\bi{d}}}$, then
\begin{align*}
\mathcal{L}_T^{(k)}\bi{v}_{\ell}&=\left(\mathcal{P}_1^{(d_1^*)}+\mathcal{P}_2^{(d_2^*-1)}
+\cdots+\mathcal{P}_{\ell-1}^{(d_{\ell-1}^*-1)}+\mathcal{Q}_{\ell}^{(d_{\ell}^*-(\ell-1))}\right)\bi{v}_{\ell}\\
&=(\ell-1)\bi{v}_{\ell}+(d_{\ell}^*-(\ell-1))\bi{v}_{\ell}=d_{\ell}^*\bi{v}_{\ell}.
\end{align*}
Now it remains to show the case that the eigenvalue is less than or equal to $d^*_{\tau_{\bi{d}}+1}$. Let $\bi{v}_{\ell}$ be the $\ell$th column of $\mathcal{L}_T^{(k)}$ where $\ell\in\{\tau_{\bi{d}}+2,\tau_{\bi{d}}+3,\cdots,k \}$. Then
$\bi{v}_{\ell}^T=[\,-\textbf{1}_{d_\ell}^T\quad\textbf{0}_{\ell-1-d_{\ell}}^T\quad d_{\ell}\quad\textbf{0}^T\,]$.
If we define $\tilde{\bi{v}}_{\ell}^T=[\,\textbf{0}_{d_\ell}^T\quad-\textbf{1}_{\ell-1-d_{\ell}}^T\quad \ell-1-d_{\ell}\quad\textbf{0}^T\,]$, then $\mathcal{L}_T^{(k)}\tilde{\bi{v}}_{\ell}=d_{\ell}\tilde{\bi{v}}_{\ell}$. Since for a threshold graph, $d_{\ell}=d_{\ell-1}^*$ for each $\ell\in\{\tau_{\bi{d}}+2,\tau_{\bi{d}}+3,\cdots,k \}$. The proof is completed.
\end{IEEEproof}
\begin{rem}\label{rem:alg}
Let $R_{ij}$ be the identity matrix with its $i$th and $j$th columns switched. If $d_i=d_j$, then $\mathcal{L}_T^{(k)}=R_{ij}\mathcal{L}_T^{(k)}R_{ij}$. Thus $\mathcal{L}_T^{(k)}\bi{v}=\lambda\bi{v}$ implies
$R_{ij}\mathcal{L}_T^{(k)}R_{ij}\bi{v}=\lambda\bi{v}$ and
\begin{equation}
\mathcal{L}_T^{(k)}R_{ij}\bi{v}=R_{ij}^{-1}\lambda\bi{v}=\lambda R_{ij}\bi{v}.
\end{equation}
Consequently, if $d_i=d_j$ and $\bi{v}$ is an eigenvector of $\mathcal{L}_T^{(k)}$, then switching the $i$th and $j$th entries of $\bi{v}$ also yields an eigenvector.
\end{rem}
\begin{lem}\label{lem:vppt}
Let
\begin{align}
\mathcal{V}_j&:=\left[\,\bi{v}_{\bar{m}^*_{j-1}+1}\,
\bi{v}_{\bar{m}^*_{j-1}+2}\,\cdots\,
\bi{v}_{\bar{m}^*_j}
\,\right]\\
\mathcal{V}_{ij}&:=\left[\,\hat{\bi{v}}_{\bar{m}^*_{j-1}+1}\,
\hat{\bi{v}}_{\bar{m}^*_{j-1}+2}\,\cdots\,
\hat{\bi{v}}_{\bar{m}^*_j}
\,\right]
\end{align}
where $\bi{v}_\ell$ is the $\ell$th column of $\mathcal{V}$ in Theorem~\ref{thm:evec} and $\hat{\bi{v}}_\ell$ is the $\bar{m}_{i-1}+1$ to $\bar{m}_i$ subvector of $\bi{v}_\ell$, $\ell\in\{\bar{m}^*_{j-1}+1,\bar{m}_{j-1}^*+2,\cdots,\bar{m}_j^*\}$. Suppose $d_{\tau_{\bi{d}}+1}=\tilde{d}_n$. Then $\mathcal{V}_{ij}$ has the following properties:
\begin{enumerate}
\item $\mathcal{V}^T_{nn}\textbf{1}=\textbf{0}\in\mathbb{R}^{(m_n-1)}$;
\item
$\mathcal{V}_{in}=\textbf{0}\in\mathbb{R}^{m_n\times (m_n-1)}$ for $i\in\left\{1,2,\cdots,\tilde{k}\right\}\setminus\{n\}$;
\item
$\mathcal{V}_{nj}=-\textbf{1}\in\mathbb{R}^{m_n\times m_j}$ for $j\in\left\{1,2,\cdots,\tilde{k}\right\}\setminus\{n\}$;

\item $\mathcal{V}_{ij}=-\textbf{1}\;\mbox{or}\;\textbf{0}
    \in\mathbb{R}^{m_i\times m_j}$, for $j\neq n$ and $i\neq j$.
\end{enumerate}
\end{lem}
\begin{IEEEproof}
These results follow directly from the definition of $\mathcal{L}_T^{(k)}$ in~(\ref{def:lkt}) and the operations of Algorithm~1 to generate $\mathcal{V}$.
\end{IEEEproof}

\begin{exmp}\label{ex:exv11}
Consider a $11$-vertex graph constructed by the string $\mathcal{S}^{(11)}=(\,0,0,1,1,0,0,0,1,0,1\,)$. Thus the degree sequence is $\bi{d}=(\,10,9,6,6,4,4,4,2,2,2,1\,)$ and the corresponding $\mathcal{L}_T^{(11)}$ is
\begin{equation*}\label{exmp}\scriptsize
-\left[\begin{array}{rrrrrrrrrrr}
-10	&1	&1	&1	&1	&1	&1	&1	&1	&1	&1\\
1	&-9	&1	&1	&1	&1	&1	&1	&1	&1	&\\
1	&1	&-6	&1	&1	&1	&1	&	&	&	&\\
1	&1	&1	&-6	&1	&1	&1	&	&	&	&\\
1	&1	&1	&1	&-4	&	&	&	&	&	&\\
1	&1	&1	&1	& &-4	&	&	&	&	&\\
1	&1	&1	&1	&	& &-4	&	&	&	&\\
1	&1	&	&	&	&&	&-2	&	&	&\\
1	&1	&	&	&	&&	&	&-2	&	&\\
1	&1	&	&	&&	&	&	&	&-2	&\\
1	&	&	&	&	&&	&	&	&	&-1
\end{array}\right].
\end{equation*}
The eigenvalues of $\mathcal{L}_T^{(11)}$ in the nonincreasing order are exactly its degree conjugate $\bi{d}^*$, namely, $(\,11,10,7,7,4,4,2,2,2,1\,)$. To derive $\mathcal{V}$ whose $i$th column corresponding to the $i$th component of $\bi{d}^*$, $\mathcal{V}$ is initialized as $\mathcal{L}^{(11)}_T$. Following step~2 to step~7 of Algorithm~1, $\mathcal{V}$ is updated and becomes
\begin{equation*}\label{exmp}\scriptsize
-\left[\begin{array}{rrrrrrrrrrr}
-10	&	&	&	&	&	&	&	&	&	&\\
1	&-8	&	&	&	&	&	&	&	&	&1\\
1	&1	&-4	&	&	&	&	&1	&1	&1	&1\\
1	&1	&1	&-3	&	&	&	&1	&1	&1	&1\\
1	&1	&1	&1	&	&1	&1	&1	&1	&1	&1\\
1	&1	&1	&1	& &-1	&1	&1	&1	&1	&1\\
1	&1	&1	&1	&	& &-2	&1	&1	&1	&1\\
1	&1	&	&	&	&&	&-5	&1	&1	&1\\
1	&1	&	&	&	&&	&	&-6	&1	&1\\
1	&1	&	&	&&	&	&	&	&-7	&1\\
1	&	&	&	&	&&	&	&	&	&-9
\end{array}\right].
\end{equation*}
Finally, after step~8 and step~9 of the algorithm we have
\begin{equation*}\scriptsize
-\left[\begin{array}{rrrrrrrrrrr}
-10	&	&	&		&	&	&	&	&	& &1\\
1	&-8	&	&		&	&	&	&	&	&1 &1\\
1	&1	&-4	&		&	&	&1	&1	&1	&1 &1\\
1	&1	&1	&-3	&	&	&1	&1	&1	&1 &1\\
1	&1	&1	&1		&1	&1	&1	&1	&1	&1 &1\\
1	&1	&1	&1	 &-1	&1	&1	&1	&1	&1 &1\\
1	&1	&1	&1		& &-2	&1	&1	&1	&1 &1\\
1	&1	&	&		&&	&-5	&1	&1	&1 &1\\
1	&1	&	&		&&	&	&-6	&1	&1 &1\\
1	&1	&	&	&	&	&	&	&-7	&1 &1\\
1	&	&	&		&&	&	&	&	&-9 &1
\end{array}\right].
\end{equation*}
\end{exmp}
Observe that the essential substrings of $\mathcal{S}^{(11)}$ are $\mathcal{S}^{(4)}=(\,0,0,1\,)$, $\mathcal{S}^{(5)}=(\,0,0,1,1\,)$ and $\mathcal{S}^{(9)}=(\,0,0,1,1,0,0,0,1\,)$, which generate $\mathbb{G}_T^{(4)}$, $\mathbb{G}_T^{(5)}$ and $\mathbb{G}_T^{(9)}$ with Laplacian matrices $\mathcal{L}_T^{(4)}$,$\mathcal{L}_T^{(5)}$ and $\mathcal{L}_T^{(9)}$ respectively. Running Algorithm~1 we obtain three sets of orthogonal eigenvectors for three matrices respectively. Namely,
 \begin{equation*} \scriptsize
-\left[\begin{array}{rrrr}
  -3 &	& 	&1 \\
1	&1	&1	&1\\
1	&-1	&1	&1\\
1	&	&-2	&1
\end{array}\right],\quad
-\left[\begin{array}{rrrrr}
 -4	&	&	&	&1\\
1	&-3	&	&	&1\\
1	&1	&1	&1	&1\\
1	&1	&-1	&1	&1\\
1	&1	&	&-2	&1\\
 \end{array}\right],
\end{equation*}
and
 \begin{equation*} \scriptsize
-\left[\begin{array}{rrrrrrrrr}
-8	&	&	&	&	&	&	&	&1    \\
1	&-4	&	&	&	&1	&1	&1	&1\\
1	&1	&-3	&	&	&1	&1	&1	&1\\
1	&1	&1	&1	&1	&1	&1	&1	&1\\
1	&1	&1	&-1	&1	&1	&1	&1	&1\\
1	&1	&1	&	&-2	&1	&1	&1	&1\\
1	&	&	&	&	&-5	&1	&1	&1\\
1	&	&	&	&	&	&-6	&1	&1\\
1	&	&	&	&	&	&	&-7	&1
 \end{array}\right].
\end{equation*}
This example illustrates that if $\bi{v}\, (\neq\textbf{1})$ is a Laplacian eigenvector of a connected threshold graph generated by the string $\mathcal{S}_1$, then an appropriate zero-padding of $\bi{v}$ must be a Laplacian eigenvector of the connected threshold graph generated by the string $\mathcal{S}_2$ which has an essential substring $\mathcal{S}_1$. We show in the next section that this property is helpful in simplifying the controllability check of a connected graph as the control matrix is binary or terminal.
\subsection{Controllability of Connected Threshold graphs}
Theorem~\ref{thm:evec} provides all independent eigenvectors required by Theorem~\ref{thm:ns} to check the controllability of $(\mathcal{L}_T^{(k)},B)$ for the general $p$-control matrix $B$. It can be readily seen that if $p$ is less than the largest multiplicity among eigenvalues of $\mathcal{L}^{(k)}_T$, then $(\mathcal{L}_T^{(k)},B)$ is uncontrollable. In the following we show that this largest multiplicity is the minimum number of controllers required to render the graph controllable and this minimum can be achieved via a binary $B$.
\begin{thm}\label{thm:ctrl} Suppose $\mathcal{S}^{(k)}$ is a constructing string that generates a $k$-vertex connected threshold graph $\mathbb{G}_T^{(k)}$ with Laplacian matrix $\mathcal{L}^{(k)}_T$, and has the essential substrings $\mathcal{S}^{(k_i)}$'s satisfying $\mathcal{S}^{(k_1)}\subseteq\mathcal{S}^{(k_2)}\subseteq\cdots\mathcal{S}^{(k_{\bar{s}})}\subseteq\mathcal{S}^{(k)}$
where $\bar{s}$ is the number of successions in $\mathcal{S}^{(k)}$ and $\mathcal{S}^{(k_i)}$ generates
$\mathbb{G}_T^{(k_i)}$ with Laplacian matrix $\mathcal{L}^{(k_i)}_T$. Suppose $B$ is bianry. $(\mathcal{L}_T^{(k)},B)$ is controllable if and only if $C^T(\mathcal{S}^{(k_i)})\hat{B}(\mathcal{S}^{(k_i)})$ is full row rank for each $i\in\mathbb{I}_{\bar{s}}$, where $C(\mathcal{S}^{(k_i)})$ and $\hat{B}(\mathcal{S}^{(k_i)})$ are local checking matrix and local control matrix, respectively, defined in Preliminaries.
\end{thm}
\begin{IEEEproof}
For each $i\in\mathbb{I}_{\bar{s}}$, Theorem~\ref{thm:evec} implies that an eigenvector of $\mathcal{L}^{(k_i)}_T$ must be a subvector of an eigenvector of $\mathcal{L}^{(k)}_T$.
If for some $i$, $C^T(\mathcal{S}^{(k_i)})\hat{B}(\mathcal{S}^{(k_i)})$ is row rank deficiency, then an eigenvector
$\bi{v}\,(\neq\textbf{1})$ of $\mathcal{L}^{(k_i)}_T$ is in the null space of $\hat{B}^T(\mathcal{S}^{(k_i)})$. Thus there exists an appropriate zero-padding vector $\tilde{\bi{v}}$ of $\bi{v}$ such that $\tilde{\bi{v}}$ is an eigenvector of $\mathcal{L}^{(k)}_T$ and is in the null space of $B^T$. Conversely, let $(\lambda,\bi{v})$ be an eigenpair of $\mathcal{L}^{(k)}_T$. If $\lambda$ is not distinct, then $\bi{v}$ must have a subvector that is also an eigenvector of  $\mathcal{L}^{(k_i)}_T$ for some $i$. Thus $\bi{v}$ in the null space of $B^T$ implies the row rank deficiency of $C^T(\mathcal{S}^{(k_i)})\hat{B}(\mathcal{S}^{(k_i)})$. If $\lambda$ is distinct, the subvector of $\bi{v}$ subject to the first two vertices can be written as either $[\,c\,c\,]^T$ or $[\,c\,-c\,]^T$, for any nonzero $c$. Suppose $\bi{v}$ is in the null space of $B^T$. In both cases, each column of $B$ subject to the first two vertices of $\mathbb{G}^{(k)}_T$ generated by $\mathcal{S}^{(k)}$ must be $[\,1\,1\,]^T$ or $[\,0\,0\,]^T$, which implies the row rank deficiency of $C^T(\mathcal{S}^{(k_1)})\hat{B}(\mathcal{S}^{(k_1)})$.
\end{IEEEproof}

\begin{cor}\label{cor:miniu}
Suppose in~Theorem~\ref{thm:ctrl} $\mathcal{L}_T^{(k)}$ has the degree sequence $\bi{d}$ in~(\ref{eq:tdseq}) and $B$ is a $p$-input control matrix. To render $(\mathcal{L}_T^{(k)},B)$ controllable,
\begin{enumerate}
\item the minimum rank of $\tilde{B}_i$, which is composed of the $(\bar{m}_{i-1}+1)$th, $\cdots$, $\bar{m}_i$th rows of $B$, is $m_i-1$, $\forall i\in\mathbb{I}_{\tilde{k}}$, for general $B$.
\item the minimum $p$ is $\max_i m_i^*$, the largest multiplicity of eigenvalues of $\mathcal{L}_T^{(k)}$ if $B$ is binary, and is $k-\tilde{k}$ if $B$ is terminal;
\end{enumerate}
\end{cor}
\begin{IEEEproof}
If the rank of $\tilde{B}_i$ is $m_i-2$ or less for some $i$, then the rank of $\tilde{C}^T(\mathcal{S}^{(k_{\sigma_i})})\tilde{B}_i$ is at most $m_i-2$. Here $\tilde{C}(\mathcal{S}^{(k_{\sigma_i})})$ is the matrix obtained by dropping the last column and the zero submatrix in the first $(m_i-1)$th columns of $C(\mathcal{S}^{(k_{\sigma_i})})$, and $\sigma_\cdot$ is the permutation function mapping $\{1,2,\cdots,\tilde{k}\}$ to $\{\sigma_1,\sigma_2,\cdots,\sigma_{\tilde{k}}\}$ for $m_i$ and $m_{\sigma_i}$ defined in Preliminaries, $i\in\mathbb{I}_{\tilde{k}}$. Consequently, $C^T(\mathcal{S}^{(k_{\sigma_i})})B(\mathcal{S}^{(k_{\sigma_i})})$ has dependent rows. Now we show how to achieve the minimum number of controllers using a binary $B$ and terminal $B$ respectively. Let $p=\max_i m_i^*$ and partition the vertices according to their degrees. We obtain the partition $\Pi:=\{\,\pi_1,\pi_2,\cdots,\pi_{\tilde{k}}\,\}$. Thus $|\pi_i|=m_i$. Let $v_{ij}$ be the $i$th vertex in $\pi_j$.
The first item of the corollary suggests that each $\pi_i$ needs at least $|\pi_i|-1$ connections to controllers. If $B$ is binary, connect $v_{ij}$ to controller $i$ for each $i\in\mathbb{I}_{m_j-2}$. If $\tilde{d}_j\neq\tau_{\bi{d}}$ then connect the $(m_j-1)$th vertex in $\pi_j$ to controller $m_j-1$, otherwise, connect the $m_j$th vertex in $\pi_j$ to controller $p$. It can be readily seen that this binary $B$ renders the graph controllable. If $B$ is terminal, each controller is connected to only one vertex, thus we need $m_i-1$ controllers for each $\pi_i$ to implement
the connection using a binary $B$. We conclude that $k-\tilde{k}$ is the minimum number of controllers for a terminal $B$.
\end{IEEEproof}

\begin{exmp} Consider the $11$-vertex threshold graph with $\mathcal{L}^{(11)}_T$ in Example~\ref{ex:exv11}. At least $3$ controllers are needed to render $(\mathcal{L}_T^{(11)},B)$ controllable. Suppose $B$ is binary. $(\mathcal{L}_T^{(11)}, B)$ is controllable if and only if $(\mathcal{L}_T^{(k)}, \hat{B}(\mathcal{S}^{(k)}))$ is controllable for each $k\in\{3,5,7\}$. Observe that $\mathcal{S}^{(k_1)}$ is always an essential substring of $\mathcal{S}^{(k_2)}$ if $k_1<k_2$. An efficient check on the controllability of $(\mathcal{L}_T^{(11)}, B)$ should work on the smaller graph first and then on the larger one, and each check requires only one rank test. In the design problem where $B$ is not available, the corollary above suggests to construct a $3$-control binary matrix $B=[\,\bi{b}_1\,\bi{b}_2\,\bi{b}_3\,]$ where $\bi{b}_1=\bi{e}_3+\bi{e}_5+\bi{e}_8$, $\bi{b}_2=\bi{e}_9$ and $\bi{b}_3=\bi{e}_7$. If $B$ is terminal, then at least $5$ controllers, e.g., $B=[\,\bi{e}_3\,\bi{e}_5\,\bi{e}_7\,\bi{e}_8\,\bi{e}_9\,]$, are needed to render $(\mathcal{L}_T^{(11)},B)$
controllable.
\end{exmp}

\subsection{Design of a Single-Input  Controllable Graph}
We now consider a special class of connected threshold graphs, namely, the class $\mathcal{C}_A$ of connected antiregular graphs. Without loss of generality, the Laplacian matrix $\mathcal{L}_A^{(k)}$ of a $k$-vertex connected antiregular graph $\mathbb{G}_A^{(k)}$ is a square matrix  of order $k$ and can be written as
\begin{equation*}{\scriptsize
-\left[\begin{array}{cccccccc}
-(k-1) &1 &\cdots  &1 &1  &\cdots &1 &1\\
1 &-(k-2) &\cdots &1  &1 &\cdots &1 &\\
\vdots &\vdots &\ddots &\vdots &\vdots &&&\\
1 &1 &\cdots &-\lfloor\frac{k}{2}\rfloor &\delta_k &&& \\
1 &1 &\cdots &\delta_k &-\lfloor\frac{k}{2}\rfloor &&&\\
\vdots &\vdots & &  &&\ddots &&\\
1 &1 & &  &&&-2&\\
1 & & & &&&&-1
\end{array}\right]}
\end{equation*}
where $\delta_k$ is $1$ if $k$ is even, and is $0$ otherwise. Observe that the repeated degree for every $k$-vertex graph in $\mathcal{C}_A$ is $\lfloor\frac{k}{2}\rfloor$.
For convenience, let
$\bar{\kappa}:=\left\lceil\frac{k}{2}\right\rceil$ and $\underline{\kappa}:=\left\lfloor\frac{k}{2}\right\rfloor$.
Clearly, we can express $k$ as $\bar{\kappa}+\underline{\kappa}$.
Combine the graphs $\mathbb{G}_A^{(\bar{\kappa})}$ and $\mathbb{G}_A^{(\underline{\kappa})}$ by adding an edge connecting the dominating vertex of $\mathbb{G}_A^{(\bar{\kappa})}$ and one of the two vertices
with the same degree in $\mathbb{G}_A^{(\underline{\kappa})}$. Call the resulting graph $\mathbb{G}^{(k)}_{\tilde{A}}$ and the corresponding Laplacian matrix $\mathcal{L}^{(k)}_{\tilde{A}}$. We have
\begin{equation}\label{def:lktd}
\mathcal{L}^{(k)}_{\tilde{A}}:=\mathcal{L}^{(k)}_c+\bi{z}_{(k)}\bi{z}_{(k)}^{T}
\end{equation}
where
\begin{equation}\label{def:Lkc}
\mathcal{L}^{(k)}_c:=\left[\begin{array}{cc}\mathcal{L}_A^{(\bar{\kappa})} &\\ &{\mathcal{L}_A^{(\underline{\kappa})}}\end{array}\right]
\end{equation}
and $\bi{z}_{(k)}:=\bi{e}_1-\bi{e}_{\bar{\kappa}+\left\lceil\frac{\underline{\kappa}}{2}\right\rceil}$
where $\bi{e}_{i}$'s are the standard basis vectors defined in Preliminaries. Occasionally we replace $\bi{z}_{(k)}$ with $\bi{z}$ for simplicity as the context is clear. Using the notations defined above, we have, for example,
\begin{equation*}
\mathcal{L}^{(8)}_{\tilde{A}}=-\left[\begin{array}{rrrrrrrr}
-4	&1	&1	&1	&	&1	&	&\\
1	&-2	&1	&	&	&	&	&\\
1	&1	&-2	&	&	&	&	&\\
1	&	&	&-1	&	&	&	&\\
	&	&	&	&-3	&1	&1	&1\\
1	&	&	&	&1	&-3	&1	&\\
	&	&	&	&1	&1	&-2	&\\
	&	&	&	&1	& 	&	&-1\\
\end{array}\right].
\end{equation*}
The following interlacing theorem follows directly from Weyl's inequality and is helpful in further spectral analysis of $\mathcal{L}^{(k)}_{\tilde{A}}$.
\begin{thm}\cite{horn13}\label{thm:intlace} Suppose $A_1$ is a Hermitian matrix of order $n$ and $A_2=A_1+\bi{z}\bi{z}^T$ where $\bi{z}$ is a nonzero column vector of size $n$. If $\lambda^{(1)}_i$ and $\lambda^{(2)}_i$ are the $i$th smallest eigenvalues of $A_1$ and $A_2$, respectively, then
\begin{align*}
\lambda_1^{(1)}\leq\lambda_1^{(2)}\leq\lambda_2^{(1)}\leq\lambda_2^{(2)}\leq\lambda_3^{(1)}\leq\cdots
\leq\lambda_n^{(1)}\leq\lambda_n^{(2)}.
\end{align*}
\end{thm}
Our objective is to combine two anti-regular graphs while maintaining the single-input controllability. To achieve this, the Laplcian matrix corresponding to the combined graph should not have any repeated eigenvalue.  In the following, we show that our proposed combination of two anti-regular graphs actually leads to distinct Laplacian eigenvalues.

\begin{lem}\label{lem:uni} The integer eigenvalues of the Laplacian matrix $\mathcal{L}^{(k)}_{\tilde{A}}$ in~(\ref{def:lktd}) are not repeated.
\end{lem}
\begin{IEEEproof} Please see Appendix for details.
\end{IEEEproof}
Define for integers $a,b$ with $a<b$,
\begin{equation}
\mathcal{R}^a_b:=R^{(-1)}_{b,b-1}\cdots R^{(-1)}_{a+2,a+1}R^{(-1)}_{a+1,a}
\end{equation}
where $R_{\alpha,\beta}^{(\gamma)}$ is the elementary row operation matrix, namely, the identity matrix except its $(\beta,\alpha)$th entry being $\gamma$. As a result,
\begin{align*}
\tilde{\mathcal{L}}^{(k)}_{\tilde{A}}(\lambda):&=
\mathcal{R}_k^{\bar{\kappa}+1}\mathcal{R}_{\bar{\kappa}}^1\left(\mathcal{L}^{(k)}_{\tilde{A}}-\lambda I\right)\\
&=\left[\begin{array}{cc}\tilde{\mathcal{L}}_A^{(\bar{\kappa})}(\lambda) &\\ &{\textbf{0}}\end{array}\right]+\left[\begin{array}{cc}\textbf{0} &\\ &\tilde{\mathcal{L}}_A^{(\underline{\kappa})}(\lambda)\end{array}\right]+\tilde{D}^{(k)}
\end{align*}
where for integer $m>2$,
\begin{equation}
\tilde{\mathcal{L}}_A^{(m)}(\lambda):=\mathcal{R}_m^1\left(\mathcal{L}_A^{(m)}-\lambda I\right)
\end{equation}
and
\begin{equation}
\tilde{D}^{(k)}:=R^{(-1)}_{\bar{\kappa}+\left\lceil\frac{\underline{\kappa}}{2}\right\rceil,
\bar{\kappa}+\left\lceil\frac{\underline{\kappa}}{2}\right\rceil-1}\bi{z}_{(k)}\bi{z}_{(k)}^T.
\end{equation}
The eigenspace property of $\mathcal{L}^{(k)}_{\tilde{A}}$, which is closely related to the controllability of $\mathbb{G}^{(k)}_{\tilde{A}}$, is given in the following theorem.
\begin{thm}\label{thm:nz} The $\left\lceil\frac{\bar{\kappa}}{2}\right\rceil$th and $\left(\left\lceil\frac{\bar{\kappa}}{2}\right\rceil+1\right)$th entries of any eigenvector of $\mathcal{L}_{\tilde{A}}^{(k)}$ in~(\ref{def:lktd}) are both nonzero.
\end{thm}
\begin{IEEEproof} Let $(\lambda,\bi{v})$ be an eigenpair of $\mathcal{L}^{(k)}_{\tilde{A}}$. The result follows trivially as $(\lambda,\bi{v})=(0,\textbf{1})$. Define
\begin{equation}\label{def:Lmda}
\Lambda:=\{0,1,\cdots,\bar{\kappa}-1\}\setminus\left\{\left\lceil\frac{\bar{\kappa}}{2}\right\rceil\right\}.
\end{equation}
If $\lambda\in\Lambda\setminus\{0\}$ and $(\lambda,\bi{v}_1)$ is an eigenpair of $\mathcal{L}^{(\bar{\kappa})}_A$, then we can write $\bi{v}^T:=\left[\,\bi{v}_1^T\;\textbf{o}^T\,\right]$ and thus the result follows from Theorem~\ref{thm:evec}. Consider the case that both $\bar{\kappa}$ and $k$ are odd.
Suppose now that $\lambda\not\in\Lambda$. By Lemma~\ref{lem:uni}, $\mathcal{L}^{(k)}_{\tilde{A}}$
has no repeated integer eigenvalue. The $\frac{\bar{\kappa}+1}{2}$th row of $\LL$ implies the equivalence of the $\left(\frac{\bar{\kappa}+1}{2}\right)$th and the $\left(\frac{\bar{\kappa}+3}{2}\right)$th entries of $\bi{v}$, and then the $\left(\frac{\bar{\kappa}-1}{2}\right)$th row of $\LL$ implies the equivalence of the $\left(\frac{\bar{\kappa}-1}{2}\right)$th, $\left(\frac{\bar{\kappa}+1}{2}\right)$th, and  $\left(\frac{\bar{\kappa}+3}{2}\right)$th entries of $\bi{v}$. Together with the $\left(\frac{\bar{\kappa}-1}{2}\right)$th row, we obtain that the $\left(\frac{\bar{\kappa}-1}{2}\right)$th, $\left(\frac{\bar{\kappa}+1}{2}\right)$th, $\left(\frac{\bar{\kappa}+3}{2}\right)$th and  $\left(\frac{\bar{\kappa}+5}{2}\right)$th entries of $\bi{v}$ are the same. Continuing these arguments
yields the expression $\left[\,\bi{v}_1^T\;\bi{v}_2^T\,\right]^T$ for $\bi{v}$ where $\bi{v}_1^T=[\,(1-\lambda)c\;c\;c\;\cdots\;c\,]\in\mathbb{R}^{\bar{\kappa}}$, for some $c\in\mathbb{R}$.
If $c=0$, then $\bi{v}_1=\textbf{0}$, and thus the first row of $\LL$ implies that the $\frac{\underline{\kappa}}{2}$th entry of $\bi{v}_2$ is $0$. In addition, the $\left(\bar{\kappa}+\frac{\underline{\kappa}}{2}\right)$th row implies that the $\left(\frac{\underline{\kappa}}{2}+1\right)$th entry is $0$ as well. Following these arguments we obtain that $\bi{v}_2=\textbf{0}$, and thus $\bi{v}=\textbf{0}$, a contradiction. We conclude that $c\neq 0$. The case of odd $k$ but even $\bar{\kappa}$, or the case of even $k$ can be proved in a similar fashion and is skipped.
 \end{IEEEproof}

 The eigenpair analysis above reveals the Laplacian controllability relation between the antiregular graph $\mathbb{G}^{(\bar{\kappa})}_A$ and the combined graph $\mathbb{G}^{(k)}_{\tilde{A}}$. We can thus design a binary control vector $\bi{b}$ that renders $\left(\mathcal{L}_{\tilde{A}}^{(k)},\bi{b}\right)$ controllable, based on the binary control vector $\bar{\bi{b}}$ that renders $\left(\mathcal{L}_{A}^{(\bar{\kappa})},\bar{\bi{b}}\right)$ controllable. A sufficient condition for the controllability of the combined graph can thus be readily derived.
 \begin{cor}\label{cor:equ} Suppose $\bar{\bi{b}}$ is a binary vector of size $\bar{\kappa}$. $\bi{b}$ is a zero-padding version of $\bar{\bi{b}}$ in the sense that $\bi{b}^T=\left[\,\bar{\bi{b}}^T\;\textbf{0}^T\,\right]$ with size $k$. Then the following statements are equivalent:
 \begin{enumerate}
 \item $\left(\mathcal{L}^{(\bar{\kappa})}_A,\bar{\bi{b}}\right)$
 is controllable;
 \item the sum of the $\left\lceil\frac{\bar{\kappa}}{2}\right\rceil$th and $\left(\left\lceil\frac{\bar{\kappa}}{2}\right\rceil+1\right)$th entries of $\bar{\bi{b}}$ is 1;
 \item $\left(\mathcal{L}^{(k)}_{\tilde{A}},\bi{b}\right)$ is controllable.
 \end{enumerate}
\end{cor}
\begin{IEEEproof}
The equivalence of the first two items follows from Theorem~\ref{thm:ctrl}. Observe that an apparent eigenpair $(\lambda,\bi{v})$ of $\mathcal{L}_{\tilde{A}}^{(k)}$ is $\lambda=\left\lceil\frac{\bar{\kappa}}{2}\right\rceil+(-1)^{\bar{\kappa}}$ and $\bi{v}=\bi{e}_{\left\lceil\frac{\bar{\kappa}}{2}\right\rceil}-\bi{e}_{\left(\left\lceil\frac{\bar{\kappa}}{2}\right\rceil+1\right)}$.
If the second item fails, then $\bi{b}$ is orthogonal to this $\bi{v}$. Suppose the second item holds. If $\lambda\in\Lambda$ in (\ref{def:Lmda}), clearly $\bi{b}$ is not orthogonal to the corresponding $\bi{v}$; otherwise, the proof of Theorem~\ref{thm:nz} suggests the form of $[\,(1-\lambda)\;1\;1\;\cdots\;1\,]$ for the first $\bar{\kappa}$ entries of $\bi{v}$ and implies the non-orthogonality of $\bi{b}$ to $\bi{v}$.
\end{IEEEproof}

\section{Conclusion}
We have studied in this work the problem of Laplacian controllability of connected threshold graphs on $k$ vertices. In the first part of the paper, we have derived a spanning set of $k$ orthogonal Laplacian eigenvectors and proposed the necessary and sufficient condition for the graphs to be controllable. If a binary control matrix is used, we have shown that we can check the rank fullness of several specific control submatrices to determine the Laplacian controllability of the graph, or use the minimum number, namely, the maximum multiplicity of entries in the conjugate of the degree sequence defining the graph, of controllers to render the graph controllable. In the second part, we have proposed to connect two antiregular graphs, one on $\lceil\frac{k}{2}\rceil$ vertices and the other on $\lfloor\frac{k}{2}\rfloor$, to form a $k$-vertex graph that is single-input controllable. The proposed graph structure is new to the class of known single-input controllable graphs and has some compromised graph parameters, compared to other well-known single-input controllable cases. For example, the diameter and the maximum vertex degree of the proposed graph are $4$ and $\lceil\frac{k}{2}\rceil$ respectively (in a path graph they are $k-1$ and $2$ respectively, and in an antiregular graph they are $2$ and $k-1$ respectively). The proposed structure enriches the class of single-input controllable graphs and provides more choices for the system designers when these parameters are their concerns. For future research topics, generalizing from the unweighted graph to the signed graph or more general graph should be interesting~\cite{sun17}. Exploring more possibilities to connect two controllable graphs while maintaining the controllability will also be the topic of interest.
\appendix[Proof of Lemma~\ref{lem:uni}]
Let $\left(\lambda_i^{(1)},\bi{v}_i^{(1)}\right)$ and $\left(\lambda_i^{(2)},\bi{v}_i^{(2)}\right)$ be the eigenpairs of $\mathcal{L}_A^{(\bar{\kappa})}$ and $\mathcal{L}_A^{(\underline{\kappa})}$ respectively, where the sequences
$\left\{\lambda_i^{(1)}\right\}_{i=1}^{\bar{\kappa}}$ and $\left\{\lambda_i^{(2)}\right\}_{i=1}^{\underline{\kappa}}$ are both nondecreasing. We can thus write $\mathcal{L}_c^{(k)}=VDV^T$ where $D=diag(\left[\,\lambda_k\;\lambda_{k-1}\;\cdots\;\lambda_1\,\right])$ and $V=[\,\bi{v}_k\,\bi{v}_{k-1}\cdots\bi{v}_1\,]$
such that
\begin{equation}
\lambda_\ell:=\left\{\begin{array}{cc}\lambda^{(1)}_{\frac{\ell+1}{2}}&\mbox{if $\ell$ is odd}\\
\lambda^{(2)}_{\frac{\ell}{2}}&\mbox{o.w.}
\end{array}\right.
\end{equation}
and
\begin{equation}
\bi{v}_\ell^T:=\left\{\begin{array}{cc}\left[\,\bi{v}^{(1)T}_{\frac{\ell+1}{2}}\;\textbf{o}^T\right]&\mbox{if $\ell$ is odd}\\
\left[\,\textbf{o}^T\;\bi{v}^{(2)T}_{\frac{\ell}{2}}\,\right]&\mbox{o.w.}
\end{array}\right..
\end{equation}
Moreover, let $\bi{x}^TA\bi{x}/\bi{x}^T\bi{x}$ be the Rayleigh-Ritz quotient for a square matrix $A$ and nonzero column vector $\bi{x}$. Equivalently, the quotient can be defined as
$R(A,\bi{x}):=\bi{x}^TA\bi{x}$ where $\|\bi{x}\|_2=1$. Suppose $\tilde{\lambda}_i$ is the $i$th smallest eigenvalue of $\mathcal{L}_{\tilde{A}}^{(k)}$. The min-max theorem has that
\begin{align}
\tilde{\lambda}_i=\min_{\dim{U}=i}\max_{\bi{x}\in U}R\left(\mathcal{L}_{\tilde{A}}^{(k)},\bi{x}\right)
=\max_{\dim{U} \atop =k-i+1}\min_{\bi{x}\in U}R\left(\mathcal{L}_{\tilde{A}}^{(k)},\bi{x}\right).
\end{align}
Thus for any $(k-i+1)$-dimensional $U$ we have $\tilde{\lambda}_i\ge\min_{\bi{x}\in U}R\left(\mathcal{L}_{\tilde{A}}^{(k)},\bi{x}\right)$, and for any $i$-dimensional space $U$ we have
$\tilde{\lambda}_i\le\max_{\bi{x}\in U}R\left(\mathcal{L}_{\tilde{A}}^{(k)},\bi{x}\right)$. If $k$ or $\bar{\kappa}$ is even, the interlacing theorem has that $\tilde{\lambda}_{2i-1}=\lambda_i$ for each $i\in\mathbb{I}_{\underline{\kappa}}$. In case $k$ and $\bar{\kappa}$ are both odd, we have
\begin{equation}
\tilde{\lambda}_{2i-1}:=\left\{\begin{array}{ll}\lambda_{2i}&\mbox{if $i\in\mathbb{I}_{\underline{\kappa}}\setminus\{\frac{\bar{\kappa}+1}{2},\}$ }\\
\lambda_{2i-1}&\mbox{if $i=\frac{\bar{\kappa}+1}{2}.$}
\end{array}\right.
\end{equation}
It remains to show that
\begin{equation}\label{ineq:tldd}
\lambda_{2i}<\tilde{\lambda}_{2i}<\lambda_{2i+1}
\end{equation}
for any $i\in\mathbb{I}_{\underline{\kappa}}$. Consider the vector $\bi{x}=V\bi{y}$ where
$\bi{y}^T=[\,y_k\;y_{k-1}\;\cdots\;y_{2i}\;\textbf{o}^T\,]$ and $\|\bi{y}\|_2=1$. Thus
\begin{align*}
R\left(\mathcal{L}_{\tilde{A}}^{(k)},\bi{x}\right)&=\bi{y}^TV^T(VDV^T+\bi{z}\bi{z}^T)V\bi{y}\\
&=\bi{y}^TD\bi{y}+\bi{y}^TV^T\bi{z}\bi{z}^TV\bi{y}\\
&=\sum_{j=2i}^{k}\lambda_j y_j^2+\left(\sum_{j=2i}^{k}y_j\bi{v}_j^T\bi{z}\right)^2.
\end{align*}
If $y_{2i+1}=y_{2i+2}=\cdots=y_k=0$, then
\begin{equation*}
R\left(\mathcal{L}^{(k)}_{\tilde{A}},\bi{x}\right)=\lambda_{2i}+(\bi{v}_{2i}^T\bi{z})^2
>\lambda_{2i},
\end{equation*}
otherwise, we have
\begin{align*}
R\left(\mathcal{L}^{(k)}_{\tilde{A}},\bi{x}\right)-\lambda_{2i}
=\sum_{j=2i+1}^k\left(\lambda_j-\lambda_{2i}\right) y_j^2+\left(\sum_{j=2i}^ky_j\bi{v}_j^T\bi{z}\right)^2>0.
\end{align*}
This shows the left strict inequality in~(\ref{ineq:tldd}). Now if we let the unit vector $\bi{y}^T=[\,\textbf{o}^T\;y_{2i}\;y_{2i-1}\;\cdots\;y_1\,]$, then
\begin{align*}
\lambda_{2i+1}-R\left(\mathcal{L}^{(k)}_{\tilde{A}},\bi{x}\right)
=&\sum_{j=1}^{2i}\left(\lambda_{2i+1}-\lambda_j\right) y_j^2-\left(\sum_{j=1}^{2i}y_j\bi{v}_j^T\bi{z}\right)^2\\
=&\;\bar{\bi{y}}^T(\bar{D}-\bar{\bi{z}}\bar{\bi{z}}^T)\bar{\bi{y}}
\end{align*}
where $\bar{\bi{y}}^T:=[\,y_1\;y_2\;\cdots\;y_{2i}\,]$ and
\begin{align*}
\bar{D}:&=\lambda_{2i+1}I-diag\left(\left[\,\lambda_1\;\lambda_2
\;\cdots\;\lambda_{2i}\,\right]\right)\,,\\
\bar{\bi{z}}^T:&=\left[\,\bi{v}_1^T\bi{z}\;\bi{v}_2^T\bi{z}\;\cdots\;\bi{v}_{2i}^T\bi{z}\,\right]=
[\,\bar{z}_1\;\bar{z}_2\;\cdots\;\bar{z}_{2i}\,].
\end{align*}
Note that the nonnegative eigenvalue of $\bar{\bi{z}}\bar{\bi{z}}^T$ is  $\bar{\bi{z}}^T\bar{\bi{z}}$ and
\begin{align*}
\bar{\bi{z}}^T\bar{\bi{z}}&\le\left[\,\bi{v}_1^T\bi{z}\;\bi{v}_2^T\bi{z}\;\cdots\;\bi{v}_k^T\bi{z}\,\right]
\left[\,\bi{v}_1^T\bi{z}\;\bi{v}_2^T\bi{z}\;\cdots\;\bi{v}_k^T\bi{z}\,\right]^T\\
   &=\bi{z}^TVV^T\bi{z}=2.
\end{align*}
By the interlacing theorem, $\bar{D}-\bar{\bi{z}}\bar{\bi{z}}^T$ is positive definite if and only if its determinant is positive. Since $\left|\bar{D}-\bar{\bi{z}}\bar{\bi{z}}^T\right|$ is equal to
\begin{equation*}
\left|\lambda_{2i+1}I-\left[\begin{array}{ccccc}
\lambda_1+\bar{z}_1^2
&\bar{z}_1\bar{z}_2 &\bar{z}_1\bar{z}_3 &\cdots &\bar{z}_1\bar{z}_{2i} \\
\bar{z}_1\bar{z}_2&\lambda_2+\bar{z}_2^2
&\bar{z}_2\bar{z}_3 &\cdots &\bar{z}_2\bar{z}_{2i} \\
\bar{z}_1\bar{z}_3 &\bar{z}_2\bar{z}_3 &\ddots &\vdots &\vdots\\
\vdots &\vdots &\cdots &\ddots &\vdots\\
\bar{z}_1\bar{z}_{2i} &\bar{z}_2\bar{z}_{2i} &\cdots  &\cdots &\lambda_{2i}+\bar{z}_{2i}^2
\end{array}\right]\right|
\end{equation*}
or
\begin{equation*}
\left|\lambda_{2i+1}I-\left[\begin{array}{ccccc}
\lambda_1+\bar{z}_1^2
&\bar{z}_1\bar{z}_2 &\bar{z}_1\bar{z}_3 &\cdots &\bar{z}_1\bar{z}_{2i} \\
\frac{\bar{z}_2}{\bar{z}_1}\left(\lambda_{2i+1}-\lambda_1\right)
&\lambda_2& & &\\
\frac{\bar{z}_3}{\bar{z}_1}\left(\lambda_{2i+1}-\lambda_1\right)
& &\lambda_3 &&\\
\vdots & & &\ddots &\\
\frac{\bar{z}_{2i}}{\bar{z}_1}\left(\lambda_{2i+1}-\lambda_1\right)
&&&&\lambda_{2i}
\end{array}\right]\right|,
\end{equation*}
we thus have
\begin{align*}
\left|\bar{D}-\bar{\bi{z}}\bar{\bi{z}}^T\right|=\left(\prod_{j=1}^{2i}
\left(\lambda_{2i+1}-\lambda_j\right)\right)
\left(1-\sum_{j=1}^{2i}\frac{\bar{z}_j^2}{\lambda_{2i+1}-\lambda_j}\right)
>0
\end{align*}
since
\begin{equation*}
\sum_{j=1}^{2i}\frac{\bar{z}_j^2}{\lambda_{2i+1}-\lambda_j}
=\left\{\begin{array}{ll}=
\bar{z}_1^2+\bar{z}_2^2=\frac{1}{\bar{\kappa}}+\frac{1}{\underline{\kappa}}<1
\quad&\mbox{if $i=1$},\\
<\sum_{j=1}^{\underline{\kappa}}\bar{z}_{2j}^2=1\quad&\mbox{if $2\leq i\leq\underline{\kappa}$}.
\end{array}\right.
\end{equation*}
This proves the right strict inequality in~(\ref{ineq:tldd}).


%



\section*{Acknowledgment}
This work is partially supported by the Ministry of Science and Technology in Taiwan under
Grant MOST-104-2221-E-005-040.


\ifCLASSOPTIONcaptionsoff
  \newpage
\fi



\bibliographystyle{IEEEtran}
\end{document}